\documentclass[11pt,oneside,french,english]{amsart}
\usepackage[T1]{fontenc}
\usepackage[latin9]{inputenc}
\usepackage[letterpaper]{geometry}
\usepackage{amsthm}
\usepackage{amssymb}
\usepackage{esint}

\makeatletter

\newcommand{\noun}[1]{\textsc{#1}}

\numberwithin{equation}{section}
\numberwithin{figure}{section}

\makeatother

\usepackage{babel}
\addto\extrasfrench{\providecommand{\fg}{\ifdim\lastskip>\z@\unskip\fi~\frqq}}

\begin{document}
\selectlanguage{french}%
\textbf{\hfill{}} \textit{Paru à }\foreignlanguage{english}{\textit{Comptes
Rendus de l'Académie des Sciences de Paris, Ser. I 347 (2009) 785\textendash{}790.}}

\selectlanguage{english}%
\vspace{0.5cm}

\begin{center}
\noun{\Large Interpolation avec contraintes sur des ensembles finis
du disque }
\par\end{center}{\Large \par}

\vspace{0.3cm}

\begin{center}
\noun{\small Rachid Zarouf}
\par\end{center}{\small \par}

\vspace{0.3cm}

\begin{flushleft}
\textbf{\small Résumé }
\par\end{flushleft}{\small \par}

{\small Etant donné un ensemble fini $\sigma$ du disque unité $\mathbb{D}=\left\{ z\in\mathbb{C}:\,\vert z\vert<1\right\} $
et une fonction $f$ holomorphe dans $\mathbb{D}$ appartenant à une
certaine classe $X$, on cherche $g$ dans une autre classe $Y$ (plus
petite que $X$) qui minimise la norme de $g$ dans $Y$ parmi toutes
les fonctions $g$ satisfaisant la condition $g_{\vert\sigma}=f_{\vert\sigma}$.
On montre que dans le cas $Y=H^{\infty}$, la constante d'interpolation
correspondante $c\left(\sigma,\, X,\, H^{\infty}\right)$ est majorée
par $a\varphi_{X}\left(1-\frac{1-r}{n}\right)$ où $n=\#\sigma$,
$r=max_{\lambda\in\sigma}\left|\lambda\right|$ et $\varphi_{X}(t)$
est la norme de la fonctionnelle d'évaluation $f\mapsto f(t)$, sur
l'espace $X.$ La majoration est exacte sur l'ensemble des $\sigma$
avec $n$ et $r$ donné.}{\small \par}

\begin{flushleft}
{\small \vspace{0.1cm}
}
\par\end{flushleft}{\small \par}

\begin{flushleft}
\textbf{\small Abstract }
\par\end{flushleft}{\small \par}

{\small Given a finite set $\sigma$ of the unit disc $\mathbb{D}=\{z\in\mathbb{C}:,\,\vert z\vert<1\}$
and a holomorphic function $f$ in $\mathbb{D}$ which belongs to
a class $X$, we are looking for a function $g$ in another class
$Y$ (smaller than $X$) which minimizes the norm $\left\Vert g\right\Vert _{Y}$
among all functions $g$ such that $g_{\vert\sigma}=f_{\vert\sigma}$.
For $Y=H^{\infty}$, and for the corresponding interpolation constant
$c\left(\sigma,\, X,\, H^{\infty}\right)$, we show that $c\left(\sigma,\, X,\, H^{\infty}\right)\leq a\varphi_{X}\left(1-\frac{1-r}{n}\right)$
where $n=\#\sigma$, $r=max_{\lambda\in\sigma}\left|\lambda\right|$
and where $\varphi_{X}(t)$ stands for the norm of the evaluation
functional $f\mapsto f(t)$ on the space $X$. The upper bound is
sharp over sets $\sigma$ with given $n$ and $r$. }{\small \par}

\begin{flushleft}
\vspace{0.1cm}

\par\end{flushleft}

\begin{flushleft}
\textbf{Abridged English version}
\par\end{flushleft}

\vspace{0.2cm}

The problem considered is the following: given $X$ and $Y$ two Banach
spaces of holomorphic functions on the unit disc $\mathbb{D}=\left\{ z\in\mathbb{C}:\,\vert z\vert<1\right\} ,$
$X\supset Y$, and a finite set $\sigma\subset\mathbb{D}$, to find
the least norm interpolation by functions of the space $Y$ for the
traces $f_{\vert\sigma}$ of functions of the space $X$, in the worst
case of $f$. 

The classical interpolation problems- those of Nevanlinna-Pick and
Carathéodory-Schur (on the one hand) and Carleson's free interpolation
(on the other hand)- are of this nature. Two first are {}``individual'',
in the sens that one looks simply to compute the norms $\left\Vert f\right\Vert _{H_{\vert\sigma}^{\infty}}$
or $\left\Vert f\right\Vert _{H^{\infty}/z^{n}H^{\infty}}$ for a
given $f$, whereas the third one is to compare the norms $\left\Vert a\right\Vert _{l^{\infty}(\sigma)}=max_{\lambda\in\sigma}\left|a_{\lambda}\right|$
and \[
inf\left(\parallel g\parallel_{\infty}:\, g(\lambda)=a_{\lambda},\:\lambda\in\sigma\right).\]

Here and everywhere below, $H^{\infty}$ stands for the space (algebra)
of bounded holomorphic functions in the unit disc $\mathbb{D}$ endowed
with the norm $\left\Vert f\right\Vert _{\infty}=sup_{z\in\mathbb{D}}\left|f(z)\right|.$
Looking at this comparison problem, say, in the form of computing/estimating
the interpolation constant\[
c\left(\sigma,\, X,\, Y\right)=sup_{f\in X,\,\parallel f\parallel_{X}\leq1}inf\left\{ \left\Vert g\right\Vert _{Y}:\, g_{\vert\sigma}=f_{\vert\sigma}\right\} ,\]
which is nothing but the norm of the embedding operator $\left(X_{\vert\sigma},\,\left\Vert .\right\Vert _{X_{\vert\sigma}}\right)\rightarrow\left(Y_{\vert\sigma},\,\left\Vert .\right\Vert _{Y_{\vert\sigma}}\right)$,
one can think, of course, on passing (after) to the limit- in the
case of an infinite sequence $\left\{ \lambda_{j}\right\} $ and its
finite sections $\left\{ \lambda_{j}\right\} _{j=1}^{n}$- in order
to obtain a Carleson type interpolation theorem $X_{\vert\sigma}=Y_{\vert\sigma}.$
But not necessarily. In particular, even the classical Pick-Nevanlinna
theorem (giving a necessary and sufficient condition on a function
$a$ for the existence of $f\in H^{\infty}$ such that $\left\Vert f\right\Vert _{\infty}\leq1$
and $f(\lambda)=a_{\lambda},$ $\lambda\in\sigma$), does not lead
immediately to Carleson's criterion for $H_{\vert\sigma}^{\infty}=l^{\infty}(\sigma).$
(Finally, a direct deduction of Carleson's theorem from Pick's result
was done by P. Koosis {[}10{]} in 1999 only). Similarly, the problem
stated for $c\left(\sigma,\, X,\, Y\right)$ is of interest in its
own. For this paper, the following question was especially stimulating
(which is a part of a more complicated question arising in an applied
situation in {[}2{]} and {[}3{]}): given a set $\sigma\subset\mathbb{D}$,
how to estimate $c\left(\sigma,\, H^{2},\, H^{\infty}\right)$ in
terms of $n=card(\sigma)$ and $max_{\lambda\in\sigma}\left|\lambda\right|=r$
only? ($H^{2}$ being the standard Hardy space of the disc).

Here, we consider the case of $H^{\infty}$ interpolation $\left(Y=H^{\infty}\right)$
and the following scales of Banach spaces $X:$

(a) $X=H^{p}=H^{p}(\mathbb{D}),$ $1\leq p\leq\infty$, the standard
Hardy spaces on the disc $\mathbb{D}$,

(b) $X=l_{a}^{2}\left(\frac{1}{(k+1)^{\alpha-1}}\right)$, $\alpha\geq1$,
the weighted spaces of all $f(z)=\sum_{k\geq0}\hat{f}(k)z^{k}$ satisfying
\[
\sum_{k\geq0}\left|\hat{f}(k)\right|^{2}\frac{1}{(k+1)^{2(\alpha-1)}}<\infty.\]

An equivalent description of this scale of spaces is:

$X=L_{a}^{2}\left(\left(1-\left|z\right|^{2}\right)^{\beta}dxdy\right)$,
$\beta=2\alpha-3>-1$, the Bergman weighted spaces of holomorphic
functions such that \[
\int_{\mathbb{D}}\left|f(z)\right|^{2}\left(1-\left|z\right|^{2}\right)^{\beta}dA<\infty.\]

For the case $\beta=0$, we shorten the notation to $X=L_{a}^{2}$.
For these two series of spaces we show \[
c_{1}\varphi_{X}\left(1-\frac{1-r}{n}\right)\leq sup\left\{ c\left(\sigma,\, X,\, H^{\infty}\right):\:\#\sigma\leq n,\,\left|\lambda\right|\leq r,\,\lambda\in\sigma\right\} \leq c_{2}\varphi_{X}\left(1-\frac{1-r}{n}\right),\]
where $\varphi_{X}(t)$, $0\leq t<1$ stands for the norm of the evaluation
functional $f\mapsto f(t)$ on the space $X$.

In order to prove the right hand side inequality, we first use a linear
interpolation: \[
f\mapsto\sum_{k=1}^{n}\left\langle f,\, e_{k}\right\rangle e_{k},\]
where $\left\langle .,.\right\rangle $ means the Cauchy sesquilinear
form $\left\langle h,\, g\right\rangle =\sum_{k\geq0}\hat{h}(k)\overline{\hat{g}(k)},$
and $\left(e_{k}\right)_{k=1}^{n}$ is the explicitly known Malmquist
basis of the space $K_{B}=H^{2}\Theta BH^{2}$, $B=\Pi_{i=1}^{n}b_{\lambda_{i}}$
being the corresponding Blaschke product, $b_{\lambda}=\frac{\lambda-z}{1-\overline{\lambda}z}$
(see N. Nikolski, {[}12{]} p. 117)). Next, we use the complex interpolation
between Banach spaces, (see H. Triebel {[}14{]} Theorem 1.9.3 p.59).
Among the technical tools used in order to find an upper bound for
$\left\Vert \sum_{k=1}^{n}\left\langle f,\, e_{k}\right\rangle e_{k}\right\Vert _{\infty}$
(in terms of $\left\Vert f\right\Vert _{X}$), the most important
is a Bernstein-type inequality $\left\Vert f^{'}\right\Vert _{p}\leq c_{p}\left\Vert B^{'}\right\Vert _{\infty}\left\Vert f\right\Vert _{p}$
for a (rational) function $f$ in the star-invariant subspace $H^{p}\cap B\overline{H}_{0}^{p}$
generated by a (finite) Blaschke product $B$, (K. Dyakonov {[}7{]}).
For $p=2$, we give an alternative proof of the Bernstein-type estimate
we need.

The lower bound problem is treated by using the {}``worst'' interpolation
$n-$tuple $\sigma=\sigma_{\lambda,\, n}=\{\lambda,\,...,\,\lambda\}$,
a one-point set of multiplicity $n$ (the Carathéodory-Schur type
interpolation). The {}``worst'' interpolation data comes from the
Dirichlet kernels $\sum_{k=0}^{n-1}z^{k}$ transplanted from the origin
to $\lambda.$ We notice that spaces $X$ of (a) and (b) satisfy the
condition $X\circ b_{\lambda}\subset X$ but this is not the case
for spaces $X$ described in (c) below for $p\neq2,$ which makes
the problem of upper/lower bound more difficult.

Other spaces considered are the following:

(c) $X=l_{a}^{p}\left(\frac{1}{(k+1)^{\alpha-1}}\right)$, $\alpha\geq1$,
$1\leq p\leq\infty$; (d) $X=L_{a}^{p}\left(\left(1-\left|z\right|^{2}\right)^{\beta}dA\right)$,
$\beta>-1$, $1\leq p\leq2.$ 

\vspace{0.1cm}

For these spaces we also found upper and lower bounds for $c\left(\sigma,\, X,\, H^{\infty}\right)$
(sometimes for special sets $\sigma$) but with some gaps between
these bounds. 

\vspace{0.2cm}

\begin{flushleft}
\textbf{A. Introduction }
\par\end{flushleft}

\vspace{0.1cm}

Le problème extremal d'interpolation est le suivant: étant donnés
$\lambda_{1},\,...,\,\lambda_{n}$ dans

\begin{flushleft}
$\mathbb{D}=\{z\in\mathbb{C}:,\,\vert z\vert<1\}$, $B=\Pi_{i=1}^{n}b_{\lambda_{i}}$
où $b_{\lambda}=\frac{\lambda-z}{1-\overline{\lambda}z}$, et $f\in\mathcal{H}ol(\mathbb{D})$,
on cherche à calculer ou estimer
\par\end{flushleft}

\[
\left\Vert f\right\Vert _{H^{\infty}/BH^{\infty}}=inf\left\{ \left\Vert g\right\Vert _{\infty}:\, f-g\in B\mathcal{H}ol(\mathbb{D})\right\} .\]

\vspace{0.2cm}

\begin{flushleft}
Les problèmes classiques de Nevanlinna-Pick (1916) et de Carathéodory-Schur
(1908), (voir {[}8{]} et {[}12{]} pour ces deux problèmes classiques
et pour des références originales), en sont des cas particuliers correspondant
respectivement à celui où les $\lambda_{j}$ sont $n$ points distincts
et au cas où $\lambda_{1}=\lambda_{2}=...=\lambda_{n}=0.$ 
\par\end{flushleft}

\begin{flushleft}
\textbf{(1) Le sujet de cette note est une version de ce dernier problème:}
\par\end{flushleft}

Trouver ou majorer/minorer la norme de la meilleure interpolante $\left\Vert f\right\Vert _{H^{\infty}/BH^{\infty}}$en
fonction de la taille de $f$ mesurée dans un espace de Banach $X$
de fonctions holomorphes dans $\mathbb{D}$, c'est à dire $\left\Vert f\right\Vert _{X}$
. De façon plus précise, il s'agit de calculer ou majorer/minorer
les constantes \[
c\left(\sigma,\, X,\, Y\right)=sup_{f\in X,\,\parallel f\parallel_{X}\leq1}inf\left\{ \left\Vert g\right\Vert _{Y}:\, g_{\vert\sigma}=f_{\vert\sigma}\right\} ,\:\mbox{et}\]

\vspace{0.03cm}
\[
C_{n,\, r}(X,Y)=sup\left\{ c(\sigma,X,Y)\,:\#\sigma\leq n\,,\,\forall j=1..n,\,\left|\lambda_{j}\right|\leq r\right\} ,\]
lorsque $Y=H^{\infty},$ $r\in[0,\,1[$ et $n\geq1.$ Nous pouvons
néanmoins faire quelques commentaires sur la constante $c\left(\sigma,\, X,\, Y\right)$
lorsque $Y$ est, comme $X,$ un espace de Banach de fonctions holomorphes
(inclus dans $X$).\vspace{0.2cm}

\begin{flushleft}
\textbf{(2) Motivations pour ce problème}
\par\end{flushleft}

\vspace{0.1cm}

\textbf{(a) }Le point de départ est une partie d'une question posée
par Laurent Baratchart (communication orale) et provenant d'un problème
d'approximation appliquée (voir {[}2{]} et {[}3{]}): trouver une estimation
de $c\left(\sigma,\, H^{2},\, H^{\infty}\right)$ en fonction de $card(\sigma)=deg(B)$
et de $max_{\lambda\in\sigma}\left|\lambda\right|=r.$

Parmi d'autres résultats, voici ci-dessous la réponse obtenue: \[
\left(\frac{\frac{1}{32}n}{1-r}\right)^{\frac{1}{2}}\leq C_{n,\, r}\left(\, H^{2},\, H^{\infty}\right)\leq\left(\frac{2n}{1-r}\right)^{\frac{1}{2}}.\]

\vspace{0.1cm}

\textbf{(b) }D'autre part, le problème de calculer/estimer $c\left(\sigma,\, X,\, Y\right)$
peut être vu comme une \textit{interpolation intermédiaire} entre
celle dite de Carleson et l'interpolation individuelle de Nevanlinna-Pick.

\vspace{0.1cm}

\textbf{(c) }On trouve une autre motivation pour l'estimation de la
constante $c\left(\sigma,\, X,\, Y\right)$ dans le calcul matriciel
où on s'interesse à la norme du calcul fonctionnel: trouver $C>0$
optimale telle que si $A$ est une matrice $n\times n$ vérifiant
$\sigma(A)\subset\sigma\subset\mathbb{D},$ telle que $\left\Vert A\right\Vert _{E\rightarrow E}\leq1$
par rapport à une certaine norme sur $\mathbb{C}^{n}$, $E=\left(\mathbb{C}^{n},\,\left|.\right|\right),$
\[
\left\Vert f(A)\right\Vert \leq C\left\Vert f\right\Vert _{\infty},\]
pour tout polynôme analytique $f$. Il est facile de voir que $C=c\left(\sigma,\, H^{\infty},\, W_{a}\right)$,
où $W_{a}$ est l'algèbre de Wiener des séries de Taylor absolument
convergentes. Notons l'apparition d'un cas intéressant pour $f\in H^{\infty}$
telle que $f_{\vert\sigma}=\frac{1}{z}_{\vert\sigma}$ (estimation
du conditionnement et des normes d'inverses des matrices $n\times n$
) ou telle que $f_{\vert\sigma}=\frac{1}{\lambda-z}_{\vert\sigma}$
(pour l' estimation de la norme de la résolvante d'une matrice $n\times n$).

Un résultat de N. Nikolski (voir {[}13{]}) nous garanti que $c\left(\sigma,\, X,\, W_{a}\right)$
est majoré par $9n.c\left(\sigma,\, X,\, H^{\infty}\right)$ pour
tout espace de Banach $X$ de fonctions holomorphes dans $\mathbb{D}$
et que cette majoration est exacte (sur $X$ et $\sigma,$ $\#\sigma\leq n$
) à une constante numérique près.

\vspace{0.2cm}

\begin{flushleft}
\textbf{(3) Les espaces considérés}
\par\end{flushleft}

\textbf{(a)} $X=H^{p}=H^{p}(\mathbb{D}),$ $1\leq p\leq\infty$, les
espaces de Hardy du disque $\mathbb{D}$,

\textbf{(b)} $X=l_{a}^{2}(\alpha):=l_{a}^{2}\left(\frac{1}{(k+1)^{\alpha-1}}\right)$,
$\alpha\geq1$, l'espace à poids des fonctions $f(z)=\sum_{k\geq0}\hat{f}(k)z^{k}$
vérifiant \[
\sum_{k\geq0}\left|\hat{f}(k)\right|^{2}\frac{1}{(k+1)^{2(\alpha-1)}}<\infty;\]
une description équivalente de cette même série d'espaces est:

$X=L_{a}^{2}\left(\beta\right):=L_{a}^{2}\left(\left(1-\left|z\right|^{2}\right)^{\beta}dxdy\right)$,
$\beta=2\alpha-3>-1$, l'espace de Bergman à poids des fonctions holomorphes
$f$ telles que \[
\int_{\mathbb{D}}\left|f(z)\right|^{2}\left(1-\left|z\right|^{2}\right)^{\beta}dxdy<\infty.\]

Pour $\beta=0$, on raccourcit la notation, $X=L_{a}^{2}$. 

\textbf{(c)} On va un peu plus loin en considérant:

$X=l_{a}^{p}(\alpha):=l_{a}^{p}\left(\frac{1}{(k+1)^{\alpha-1}}\right)$,
$\alpha\geq1$, $1\leq p\leq\infty$, puis $X=L_{a}^{p}\left(\beta\right):=L_{a}^{p}\left(\left(1-\left|z\right|^{2}\right)^{\beta}dxdy\right)$,
$\beta>-1$, $1\leq p\leq2.$

\vspace{0.3cm}

\begin{flushleft}
\textbf{B. Ce qui est montré}
\par\end{flushleft}

\vspace{0.1cm}

Nous commençons par étudier le cas d'espaces de Banach généraux $X$
et $Y$ verifiant les propriétés naturelles suivantes:\def\theequation{$P_{1}$}\begin{equation}
pour\: tout\:\epsilon>0,\: Hol((1+\epsilon)\mathbb{D})\: est\: continument\: inclus\: dans\: Y\label{eq:}\end{equation}

\def\theequation{$P_{2}$}

\begin{flushleft}
\begin{equation}
Pol_{+}\subset X\: et\; Pol_{+}\: est\: dense\: dans\: X,\label{eq:}\end{equation}
où $Pol_{+}$ désigne l'espace des polynômes analytiques à coefficients
complexes $p$, $p(z)=\sum_{k=0}^{N}a_{k}z^{k},$
\par\end{flushleft}

\def\theequation{$P_{3}$}

\begin{flushleft}
\begin{equation}
\left[f\in X\right]\Rightarrow\left[z^{n}f\in X\,,\,\forall n\geq0\: and\:\overline{lim}\left\Vert z^{n}f\right\Vert ^{\frac{1}{n}}\leq1\right],\label{eq:}\end{equation}

\par\end{flushleft}

\def\theequation{$P_{4}$}\begin{equation}
\left[f\in X,\,\lambda\in\mathbb{D},\, and\, f(\lambda)=0\right]\Rightarrow\left[\frac{f}{z-\lambda}\in X\right].\label{eq:}\end{equation}

\begin{flushleft}
\textbf{Lemme 0.}\textit{ Soient $X,Y$ deux espaces de Banach vérifiant
les propriétés $\left(P_{i}\right)$, $i=1...4$ . Pour tout $n\geq1$,
$r\in[0,\,1)$, on a $C_{n,\, r}(X,Y)<\infty.$}
\par\end{flushleft}

Puis, en étudiant le cas particulier où $Y=H^{\infty}$ et où $X$
parcourt les espaces décrits en A.(3)-(a) et A.(3)-(b), nous obtenons
des majorations/minorations du type

\[
c_{1}\varphi_{X}\left(1-\frac{1-r}{n}\right)\leq C_{n,\, r}\left(X,\, H^{\infty}\right)\leq c_{2}\varphi_{X}\left(1-\frac{1-r}{n}\right),\]
 où $\varphi_{X}(t)$, $0\leq t<1$ est la norme de la fonctionnelle
d'évaluation $f\mapsto f(t)$ sur l'espace $X$.

\vspace{0.2cm}

Plus précisément, en ce qui concerne les espaces de Hardy $H^{p}$
du A.(3)-(a), on obtient le théorème suivant.

\begin{flushleft}
\textbf{Théorème 1. }\textit{Soit $p\in2\mathbb{Z}_{+}$. Il existe
une constante $A_{p}$ dépendant de $p$ uniquement telle que pour
tout $n\geq1$, $r\in[0,\,1),$ \[
\frac{1}{32^{\frac{1}{p}}}\left(\frac{n}{1-r}\right)^{\frac{1}{p}}\leq C_{n,\, r}\left(H^{p},H^{\infty}\right)\leq A_{p}\left(\frac{n}{1-r}\right)^{\frac{1}{p}}.\]
De plus, la majoration est vraie pour tout réel $p,$ $1\leq p\leq+\infty$.}
\par\end{flushleft}

\vspace{0.2cm}
Quant au cas des espaces à poids (où de façon équivalente celui des
espaces de Bergman à poids radial) du A.(3)-(b), $X=l_{a}^{2}\left(\alpha\right)=L_{a}^{2}\left(2\alpha-3\right),$
on obtient l'encadrement suivant.

\begin{flushleft}
\textbf{Théorème 2.}\textit{ Soit $\alpha\geq1$ tel que $2\alpha-1$
soit entier. Il existe des constantes $a$ et $A$ telles que pour
tout $n\geq1,r\in[0,\,1),$}
\par\end{flushleft}

\textit{\[
a\left(\frac{n}{1-r}\right)^{\frac{2\alpha-1}{2}}\leq C_{n,\, r}\left(l_{a}^{2}\left(\alpha\right),\, H^{\infty}\right)\leq A\left(\frac{n}{1-r}\right)^{\frac{2\alpha-1}{2}},\]
où les constantes $a$ et $A$ sont telles que $a\asymp\frac{1}{2^{3N}(2N)!}$
et $A\asymp N^{2N}$, $N$ étant la partie entière de $\alpha$. De
plus, la majoration est vraie pour tout réel $\alpha\geq1.$ (La notation
$x\asymp y$ signifie qu'il existe des constantes numériques $c_{1},\, c_{2}>0$
telles que $c_{1}y\leq x\leq c_{2}y$). }

\vspace{0.2cm}

Enfin, en ce qui concerne le cas des espaces $X$ du A.(3)-(c), les
résultats obtenus sont plus faibles et ne répondent pas, comme c'était
le cas précédemment, à la conjecture faisant intervenir $\varphi_{X}$
définie ci-dessus. Nous donnons néanmoins des majorations/minorations
pour la quantité $c\left(\sigma,\, X,\, H^{\infty}\right),$ parfois
pour des $\sigma$ spécifiques et avec des écarts entre les bornes
intervenant dans ces encadrements.

\begin{flushleft}
\textbf{Théorème 3.} \textit{(1) Soit $\alpha\geq1$ et $X=l_{a}^{p}\left(\alpha\right)$$.$
Pour $1\leq p\leq+\infty$ , il existe une minoration de $C_{n,\, r}\left(X,\, H^{\infty}\right)$
de l'ordre de $1/(1-r)^{\alpha-1/p}$. Pour $1\leq p\leq2$ (resp.
$2\leq p\leq+\infty$ ), il existe une majoration de $C_{n,\, r}\left(X,\, H^{\infty}\right)$
de l'ordre de $\left(\frac{n}{1-r}\right)^{\alpha-1/2}$ (resp. de
l'ordre de $\left(\frac{n}{1-r}\right)^{\alpha+1/2-2/p}$). }
\par\end{flushleft}

\textit{(2) Soient $\lambda\in\mathbb{D}$, $\beta>-1$, $1\leq p\leq2$
et $X=L_{a}^{p}\left(\beta\right)$. Alors il existe une majoration
de $c\left(\sigma_{\lambda,\, n},\, X,\, H^{\infty}\right)$ de l'ordre
de $\left(\frac{n}{1-\left|\lambda\right|}\right)^{(\beta+2)/p}.$}

\vspace{0.3cm}

\begin{flushleft}
\textbf{C. Les moyens pour montrer cela}
\par\end{flushleft}

\vspace{0.1cm}

\textbf{(a) Majorations}

\vspace{0.1cm}

\textbf{(i)} Nous avons choisi d'utiliser une interpolation linéaire\[
T:\, f\mapsto\sum_{k=1}^{n}\left\langle f,\, e_{k}\right\rangle e_{k},\]
où $\left\langle .,.\right\rangle $ est la forme sesquilinéaire de
Cauchy $\left\langle h,\, g\right\rangle =\sum_{k\geq0}\hat{h}(k)\overline{\hat{g}(k)},$
et $\left(e_{k}\right)_{k=1}^{n}$ est la base de l'espace $K_{B}=H^{2}\Theta BH^{2}$,
dite de Malmquist, connue de façon explicite, (voir N. Nikolski, {[}12{]}
p. 117)). Ce choix est justifié par le fait que si $X=H^{2}$, l'opérateur
$T$ coïncide avec la projection orthogonale de $H^{2}$ sur $H^{2}\Theta BH^{2}.$
Nous conserverons ce choix même dans le cas plus général où $X=H$
est un espace de Hilbert différent de $H^{2}$, car pour ce type d'espace
la projection orthogonale de $H$ sur $H\Theta BH$ demeure implicite.
D'autre part, si $X$ n'est pas un espace de Hilbert, trouver la {}``meilleure''
interpolante de $f$ est encore moins clair, il y aura donc un prix
à payer relativement à notre choix. A ce propos, en général la vraie
interpolation optimale est non-linéaire, voir S.A. Vinogradov ({[}15{]})
pour les détails.

\vspace{0.1cm}

\textbf{(ii)} Pour le cas A.(3)-(a) (Théorème 1) avec $p=2$, on profite
du fait que $e_{k}\in\mathcal{H}ol(\left|z\right|<\frac{1}{r})$ où
$r=max_{\lambda\in\sigma}\left|\lambda\right|$ pour majorer $\left\Vert g\right\Vert _{\infty}.$
Pour généraliser au cas $p\geq1$ quelconque on utilise un résultat
d'interpolation de P. Jones notamment que $[H^{1},\, H^{\infty}]_{\theta}=H^{p},$
voir {[}9{]}.

\vspace{0.1cm}

Pour le cas A.(3)-(b) (Théorème 2), on utilise les points ci-après.

\textbf{(iii) }Pour $X=l_{a}^{2}\left(N+1\right)$, on fait apparaître
la quantité $\left\Vert g^{(N)}\right\Vert _{H^{2}}$ (qui est comparable
à $\left\Vert g^{}\right\Vert _{X^{\star}}$) que l'on majore en fonction
de $\left\Vert f\right\Vert _{H^{2}}$ à l'aide d'une inégalité type
Bernstein sur les fonctions rationelles à pôles dans $^{c}\overline{\mathbb{D}}$,
que l'on montrera. Plus précisément, on démontrera le lemme suivant,
qui est un analogue pour le disque d'un résultat de K. Dyakonov démontré
dans le demi-plan, voir {[}7{]}.

\begin{flushleft}
\textbf{Lemme 1.}\textit{ Soit $g\in K_{B}=:H^{2}\Theta BH^{2}.$
Alors\[
\left\Vert g^{'}\right\Vert _{H^{2}}\leq\frac{5}{2}\frac{n}{1-r}\left\Vert g\right\Vert _{H^{2}},\]
où comme toujours, $r=max_{\lambda\in\sigma}\left|\lambda\right|.$
Par récurrence,\[
\left\Vert g^{(k)}\right\Vert _{H^{2}}\leq k!\left(\frac{5}{2}\right)^{k}\left(\frac{n}{1-r}\right)^{k}\left\Vert g\right\Vert _{H^{2}},\]
pour tout $k=0,\,1,\,...$ .}
\par\end{flushleft}

Notre preuve est différente de celle de M. Dyakonov et elle donne
en particulier une constante ($5/2$) plus petite. En général, on
notera que les inégalités type Bernstein ont déjà fait l'objet de
nombreuses publications. Entre autres, le chapitre 7 du livre de P.
Borwein et T. Erdélyi, voir {[}5{]}, y est consacré.\textit{ }C'est
aussi le cas de la thèse de A. Baranov, voir {[}1{]}, et de l'ouvrage
de R. A. DeVore and G. G. Lorentz, voir {[}6{]}.

\vspace{0.1cm}

\textbf{(iv) }Enfin, comme en (ii), on interpole entre $l_{a}^{2}\left(N\right)$
et $l_{a}^{2}\left(N+1\right)$ (interpolation classique complexe
entre espaces de Banach, voir {[}14{]} ou {[}4{]}).

\vspace{0.1cm}

\textbf{(v) }Pour traiter le cas A.(3)-(c), on fait de même qu'en
(iv) mais entre $l_{a}^{p_{1}}\left(\alpha\right)$ et $l_{a}^{p_{2}}\left(\alpha\right)$,
puis entre $L_{a}^{p_{1}}\left(\beta\right)$ et $L_{a}^{p_{2}}\left(\beta\right).$

\vspace{0.3cm}

\textbf{(b) Minorations}

\vspace{0.1cm}

On se rend compte, grâce l'interpolation de Carleson, que la séquence
la ''pire'' est probablement $\sigma=\sigma_{\lambda,\, n}=\{\lambda,\,...,\,\lambda\}$
($n$ fois). En effet, dans ce cas la constante de Carleson explose
et, au travers de la majoration, (qui est vraie pour toute séquence
$\sigma$ de $\mathbb{D}$),

\[
c(\sigma,\, X,\, H^{\infty})\leq C_{I}(\sigma).max_{1\leq i\leq n}\left\Vert \varphi_{\lambda_{i}}\right\Vert ,\]
où $\varphi_{\lambda}(f)=f(\lambda)$ et \[
C_{I}(\sigma)=sup_{\parallel a\parallel_{l^{\infty}}\leq1}inf\left(\parallel g\parallel_{\infty}:\, g\in H^{\infty},\, g_{\vert\sigma}=a\right),\]
est la constante de Carleson relative à $\sigma,$ on comprend que
pour $\sigma$ ayant une constante d'interpolation $C_{I}(\sigma)$
{}``raisonnable'', la quantité $c(\sigma,\, X,\, H^{\infty})$ se
comporte comme $max_{i}\left\Vert \varphi_{\lambda_{i}}\right\Vert $.
En revanche, pour des sequences $\sigma$ {}``serrées'', la constante
$C_{I}(\sigma)$ peut être si grande que la dernière majoration peut
en devenir très grossière . 

\vspace{0.1cm}

\textbf{(i)} On remarque d'abord que 

\begin{flushleft}
 \[
c\left(\sigma_{0,\, n},\, H^{2},\, H^{\infty}\right)\geq\frac{1}{\left\Vert p_{n}\right\Vert _{H^{2}}}\left\Vert p_{n}\star K_{n}\right\Vert _{\infty}\geq\frac{1}{\sqrt{n-1}}\left(p_{n}\star K_{n}\right)(1)\geq\frac{1}{2}\sqrt{n},\]
 où $K_{n}$ désigne le noyau de Fejer d'ordre $n,$ et $p_{n}=\sum_{k=0}^{n-1}z^{k}$.
\par\end{flushleft}

\textbf{(ii)} Le cas $c\left(\sigma_{0,\, n},\, X,\, H^{\infty}\right)$
pour $X$ espace de Hilbert du A.(3)-(b) se traite de la même façon,
en remplaçant $p_{n}$ par une puissance de $p_{n}.$

\textbf{(iii) }La minoration de $c(\sigma_{\lambda,\, n},\, H^{2},\, H^{\infty})$
se {}``déduit'' de (ii) en considérant, au lieu de la fonction $p_{n}=\sum_{k=0}^{n-1}z^{k}$,
une fonction construite à partir de $p_{n}\circ b_{\lambda}$, ou
plus simplement $p_{n}\circ b_{r},$ $0\leq r<1$ puisque les normes
considérées sont invariantes par rotation. On prouve ainsi la minoration
du Théorème 1 pour $p=2$.

\textbf{(iv)} Encore une fois, le cas $c(\sigma_{\lambda,\, n},\, X,\, H^{\infty})$
pour $X$ espace de Hilbert du A.(3)-(b) se traite comme en (iii),
en remplaçant $p_{n}$ par une puissance de $p_{n}$. L'observation
principale réside dans le fait que $X=l_{a}^{2}\left(\alpha\right)=\varphi\left(H^{2}\right)=H(\varphi\circ K)$\noun{\Large{}
}dans le sens d'Aronszajn-deBranges, (voir {[}12{]} p.320, point $(k)$
de l'Exercice 6.5.2), avec $\varphi(z)=z^{2\alpha-1}$ et $K(\lambda,z)=k_{\lambda}(z)=\frac{1}{1-\bar{\lambda}z}$.
En particulier, on utilise l'inégalité suivante, vraie pour tout $f\in H^{2}$:

\[
\left\Vert \varphi\circ f\right\Vert {}_{X}^{2}\leq\varphi\left(\left\Vert f\right\Vert {}_{2}^{2}\right).\]

Il est bon de noter que l'opérateur de composition par $b_{\lambda}$
stabilise les espaces du A.(3)-(a) et du A.(3)-(b) mais qu'en revanche
ce n'est pas le cas pour les espaces $l_{a}^{p}\left(\alpha\right)$
pour $p\neq2$ ce qui rend le problème de minoration mais aussi de
majoration plus difficile.

\textbf{(v) }Enfin, pour montrer la minoration du Théorème 3, on utilise
simplement que \[
C_{n,\, r}(X,H^{\infty})\geq\left\Vert \psi_{r}\right\Vert _{X},\]
où $\psi_{r}(f)=f(r).$ (Mais la minoration ainsi obtenue n'est pas
optimale).

\vspace{0.1cm}

\begin{flushleft}
\textbf{Remerciements}
\par\end{flushleft}

\vspace{0.1cm}

Je tiens à remercier chaleureusement le Professeur Nikolai Nikolski
pour son aide inestimable.

\[
\begin{array}{c}
\mbox{Equipe\:\ d'Analyse\:\ et\:\ Géométrie,}\\
\mbox{Institut\:\ de\:\ Mathématiques\:\ de\:\ Bordeaux,}\\
\mbox{\mbox{Université\:\ Bordeaux,\:351\:\ Cours\:\ de\:\ la\:\ Libération,\,33405\:\ Talence,\:\ France}.}\\
\mbox{E-mail\:\ address:\:\ rzarouf@math.u-bordeaux1.fr}\end{array}\]

\end{document}